# A VOLUME PRODUCT REPRESENTATION AND ITS RAMIFICATIONS IN $l_p^n, 1 \le p \le \infty$

## D. KARAYANNAKIS

**Abstract**

If $|B_p^n|, 1 < p < \infty$, is the volume of the unit $p$-ball in $R^n$ and $q$ the Hölder conjugate exponent of $p$, we represent the product $|B_p^n||B_q^n|$ as a suitable function free of its gamma symbolism; this representation will allows us to confirm by use of basic classical analysis tools, for the particular case $K = B_p^n, 1 \le p \le \infty$, two conjectured and/or proved lower and upper bounds for the volume product of centrally symmetric convex bodies of the Euclidean $R^n$ that play a central role in convex geometric analysis.

## Part I: Notation, goals and main strategy

One of the key notions in convex analysis is the volume product $M(K) = |K||K^o|$ where $K$ is a (centrally symmetric) convex compact set of the Euclidean $(R^n, \langle .,. \rangle)$ with non-empty interior (or simply "body") and $K^o = \{y \in R^n : |\langle x, y \rangle| \le 1\}, \forall x \in K\}$ is the so-called polar set of $K$

For this product, one of the long standing conjectures stated by Mahler (who proved it for $n = 2$), claims that $\dfrac{4^n}{n!} \le M(K)$ for origin symmetric bodies $K$. This conjecture has been confirmed in many special cases of $K$ and in particular for bodies symmetric with respect to the coordinate planes, which naturally include the $p$-balls, using advanced Banach space theory. On the other hand, for an upper bound, we have the inequality $M(K) \le M(B_2^n)$ proved in 1948 by Santaló (and much earlier for $n = 2$ by Blashke). A survey of the above facts and other related results can be found e.g. in [4].

We set forward in this short work to establish Mahler's conjecture and the Blashke-Santaló inequality for the case $K = B_p^n, 1 \le p \le \infty$, by using only basic special functions and classical analysis theory.

It is clear that $(B_p^n)^o = B_q^n$, where $\dfrac{1}{p} + \dfrac{1}{q} = 1$, and also directly verifiable through multiple integration (see e.g. [4] and [1]) that $|B_p^n| = \dfrac{2^n \Gamma(1+\frac{1}{p})^n}{\Gamma(1+\frac{n}{p})}$ for $1 < p < \infty$ and $|B_1^n| = |B_\infty^n| = 2^n$.

So we should be able to manipulate the expression $M(B_p^n) = 4^n \dfrac{[\Gamma(1+\frac{1}{p})\Gamma(1+\frac{1}{q})]^n}{\Gamma(1+\frac{n}{p})\Gamma(1+\frac{n}{q})}$. (1)

Thus, in Part II, we start with a suitable for our goals gamma functions ratio result (Lemma II.1) that will allow us to represent $M(B_p^n)$ as a suitable function $M(n,p)$, for $1 < p < \infty$, free of its gamma symbolism(Prop. II.1); in Part III, by establishing $\frac{dM(n,p)}{dp} \geq 0$ for $1 < p \leq 2$ (Lemma.III.1) and then using the evident facts that $M(n,p) = M(n,q)$ and $1 < q \leq 2$ iff $2 \leq p < \infty$, and also by examining separately the case $p=1$, we will obtain the announced results (Prop.III.1).We also obtain, as byproducts of independent interest, two seemingly new closed formulae concerning infinite products(Corol.III.1 & Corol.III.2).

## Part II: A Lemma and a Proposition

Lemma II.1

For $x > 0$ and $a < 1$ we have that $\Gamma(1-a)\frac{\Gamma(x+a)}{\Gamma(x)} = \prod_{k=1}^{\infty} \frac{k(k+x-1)}{(k-a)(k+x+a-1)}$.

Proof :
Let $P(x,a)$ the above infinite product .One of the various classical definitions for the gamma function is $\Gamma(z) = \lim_{k \to \infty} \frac{k^{z-1}k!}{(z)_k}$, for any complex $z \neq 0, -1, -2, \ldots$, where by $(z)_k$ we have denoted the shifted factorial $z(z+1)\ldots(z+k-1)$ .Substituting in the above limit respectively $z = 1-a, z = x+a,$ and $z = x$ and after simplifying we arrive at

$\lim_{k \to \infty} \frac{k!(x)_k}{(1-a)_k (x+a)_k}$ which is identical to $P(x,a)$. □

Remark II.1
(i)Independently of the above argumentation we can easily check that $P(x,a)$ exists as a two-variable function over $(0, \infty) \times [0,1)$ since evidently $0 \leq P(x,a)$ and, by use of the inequality $\log t \leq t - 1, t > 0$, $\log P(x,a) \leq a(x+a-1) \sum_{k=1}^{\infty} \frac{1}{(k-a)(k+x+a-1)} < \infty$.
(ii)Lemma II.1 was proved (formally and from "scratch") in [3] where the scheme of the proof served different purposes concerning the numerical evaluation of the gamma and psi function).

**Proposition II. 1**

For $1 < p < \infty$, $M(n,p) = 4^n h(p)^{2n-2} \prod_{k=1}^{\infty} \frac{(k^2 + k + \frac{1}{pq})^{n-2} (k^2 + nk + \frac{n^2}{pq})}{(k+1)^{2n-2}}$ (2)

where $h(p) = \frac{\pi}{pq \sin(\frac{\pi}{p})}$.

Proof :

Evidently $h(1^+) = 1$ and so (2) is trivially true for $n = 1$.Thus we can consider $n \geq 2$.

At first let us rewrite (1) as $M(n,p) = \frac{4^n}{n^2} \frac{1}{(pq)^{n-1}} \frac{[\Gamma(\frac{1}{p})\Gamma(\frac{1}{q})]^n}{\Gamma(\frac{n}{p})\Gamma(\frac{n}{q})}$.

Based on Lemma II.1 for $a = \frac{1}{p}$ and $x = \frac{m}{p}, m = 1,,...,n-1$ we see that

$$\Gamma(\frac{n}{p}) = \frac{\Gamma(\frac{1}{p})}{\Gamma(\frac{1}{q})^{n-1}} \prod_{m=1}^{n-1} P(\frac{n-m}{p}, \frac{1}{p}) \tag{3}$$

where $P(x,a)$ was defined above in Lemma II.1. Working in a similar way we obtain the conjugate expression for $\Gamma(\frac{n}{q})$. We observe now that

$$P(\frac{n-m}{p}, \frac{1}{p}) \, P(\frac{n-m}{q}, \frac{1}{q}) = \frac{1}{pq}(\frac{n-m}{n+1-m})^2 \prod_{k=1}^{\infty} \frac{(k+1)^2}{g_k(p)} \frac{r_k(n-m, p)}{r_k(n+1-m, p)} \tag{4}$$

where we have set $r_k(j, p) = k^2 + jk + \frac{j^2}{pq}$ and $g_k(p) = r_k(1, p)$.

Combining (1), (3), and (4), and after "telescoping" inside the infinite product, we obtain the announced result. □

Part III: One more Lemma and a conclusive Proposition

**Lemma III.1**

Let $1 < p \leq 2$. Then $\frac{dM(n, p)}{dp} \geq 0$ for all $n$ with equality iff $p = 2$ for $n \neq 1$.

Proof:

$M(1, p) = 4^n$ and so once more we can consider $n \geq 2$. By straightforward differentiation with respect to $p$ and in the case of the infinite product of (2) by logarithmic differentiation (noticing that this product by construction is a real analytic function of $p$ having as logarithm a uniformly converging series of differentiable functions of $p$) we obtain

$$\text{sgn}[\frac{dM(n, p)}{dp}] = \text{sgn}[\,(2n-2)h' + (\frac{1}{pq})'h \sum_{k=1}^{\infty}(\frac{n-2}{g_k} + \frac{n^2}{r_k})]. \tag{5}$$

In (5) $h = h(p)$ as defined in (2) and $g_k = g_k(p), r_k = r_k(n, p)$ as defined in (4); still for simplicity from now on we suppress any possible dependence from $p$ and $n$.

Since evidently the infinite sum in (5) and $h$ are positive and $(pq)' = \frac{p(p-2)}{(p-1)^2} \leq 0$ (with equality iff $p = 2$) it will suffice to show that $h' \geq 0$ (with equality iff $p = 2$).

Simple calculations show that $\text{sgn}\, h' = -\text{sgn}[pq \sin(\frac{\pi}{p})]' = \text{sgn}(\omega)$, where

$$\omega = \omega(p) = p(2-p)\sin(\frac{\pi}{p}) + \pi(p-1)\cos(\frac{\pi}{p}). \tag{6}$$

Now $\text{sgn}(\omega') = \text{sgn}[\pi(1-p)\sin(\frac{\pi}{p}) + (p^2 + 2p - 2)\cos(\frac{\pi}{p})] < 0$, \tag{7}

and since $\omega(2) = 0$ we conclude that $\omega \geq 0$ with equality iff $p = 2$. □

We are now ready to prove the announced double inequality concerning $M(n, p)$.

**Proposition III.1**

For $1 \leq p \leq \infty$ and all positive integers $n$, $\dfrac{4^n}{n!} \leq M(n,p) \leq M(n,2)$, with left (resp. right) equality iff $p = 1$ (resp. $p = 2$), whenever $n \neq 1$.

Proof:
Based on Lemma III.1 we can see immediately that for $1 < t < p < \infty$,
$M(n,t) < M(n,p) \leq M(n,2)$, with equality iff $p = 2$.
We only need to examine (2) whenever $p \to 1^+$. It is clear that $M(n,1^+) = 4^n s_n$, with

$$s_n = \prod_{k=1}^{\infty} \frac{k^{n-1}}{(k+1)^n}(k+n). \text{ Since } s_1 = 1 \text{ and noticing that } \frac{s_{n+1}}{s_n} = \prod_{k=1}^{\infty} \left(\frac{k}{k+1}\right)\left(\frac{k+n+1}{k+n}\right) = \frac{1}{n+1}$$

we see that $s_n = \dfrac{1}{n!}$ and we are done. □

**Remark III.1**

Independently of the above approach one can easily verify that $s_n$ converges for all $n$, by the logarithmic series test along the lines of Remark II.1(i): evidently $s_n > 0$ and also for $n \geq 2$ (by the crudest logarithmic inequality) $\log s_n < -\sigma_n$, where $\sigma_n = \sum_{m=2}^{n} \binom{n}{m} \sum_{k=1}^{\infty} \dfrac{k^{n-m}}{(k+1)^n} < \infty$. Thus we obtain the (crude but not obvious) inequality

$\dfrac{1}{n!} < e^{-\sigma_n}$ for $n \geq 2$.

We conclude this work with two closed type formulae that are automatically true when we set $p = 2$ in (2). Exploiting the fact that $M(n,2) = \left|B_2^n\right|^2 = \dfrac{4\pi^n}{n^2 \Gamma^2\left(\dfrac{n}{2}\right)}$ and the classical value of

$\Gamma\left(\dfrac{n}{2}\right)$ for even and odd $n$ (see e.g.[2]) we arrive, respectively, for $m \in N$ at the following:

**Corollary III.1**

$$\prod_{k=1}^{\infty} \frac{(2k+1)^{2m-2}(2k+2m)}{(2k+2)^{2m-1}} = \frac{\left(\dfrac{4}{\pi}\right)^{m-1}}{m!}.$$

Corollary III.2

$$\prod_{k=1}^{\infty} \frac{(2k+1)^{2m-1}(2k+2m+3)}{(2k+2)^{2m}} = \frac{\left(\dfrac{8}{\pi}\right)^{m+1}}{(2m+3)!!}.$$

**Remark III.2**
It is clear that the above two formulae provide us in closed form the value of the seemingly untabulated infinite product $\prod_{k=1}^{\infty} \dfrac{(2k+2)(2k+2m)}{(2k+1)(2k+2m+3)}$.

DEPARTMENT OF SCIENCES ,TEI OF CRETE, 71004 HERAKLION, GREECE
*E-mail address*: dkar@stef.teiher.gr